# A STUDY ON DUAL MANNHEIM PARTNER CURVES

M. A. GUNGOR, M. TOSUN

*Sakarya University, Faculty of Arts and Science*

*Department of Mathematics, Sakarya / TÜRKİYE*

**Abstract:** Mannheim partner curves are studied by Liu and Wang [1,2]. Orbay and others extended the theory of the Mannheim curves to the ruled surface in Euclidean 3-space $E^3$ [3]. We obtain the relationships between the curvatures and the torsions of the dual Mannheim partner curves with respect to each other.

## 1- Introduction

Ruled surfaces are the surfaces which are generated by moving a straight line continuously in the space and are one of the most important topics of differential geometry. A ruled surface can always be described (at least locally) as the set of points swept by a moving straight line. These kinds of surfaces are used in many areas of sciences such as Computer Aided Geometric Design (CAGD), mathematical physics, moving geometry and kinematics for modeling the problems. So, many geometers and engineers have studied on ruled surfaces in different spaces such as Euclidean space [4-7]. They have investigated and obtain many properties of the ruled surfaces in this space.

Since building materials such as wood are straight, they can be considered as straight lines. The result is that if engineers are planning to construct something with curvature, they can use a ruled surface since all the lines are straight [3].

The curves are a fundamental structure of differential geometry. In the differential geometry of a regular curve in the Euclidean 3-space $E^3$, it is well-known that one of the important problem is the characterization of a regular curve. The curvature functions $k_1$ (curvature $\chi$) and $k_2$ (torsion $\tau$) of a regular curve play an important role to determine the shape and size of the curve [7-9]. For example: If $k_1 = k_2 = 0$, then the curve is a geodesic. If $k_1 \neq 0$ (constant) and $k_2 = 0$, then the curve is a circle with radius $1/k_1$. If $k_1 \neq 0$ (constant) and $k_2 \neq 0$ (constant), then the curve is a helix in the space, etc.

Another way to classification and characterization of curves is the relationship between the Frenet vectors of the curves. For example, in the plane, a curve $\alpha$ rolls on a straight line, the center of curvature of its point of contact describes a curve $\beta$ which is the Mannheim of $\alpha$ Mannheim partner curves in three dimensional Euclidean 3-space are studied by Liu and Wang [1,2]. They have given the definition of Mannheim offsets as follows: Let $C$ and $C_1$ be two space curves. $C$ is said to be a Mannheim partner curve of $C_1$ if there exists a one to one correspondence between their points such that the binormal vector of $C$ is the principal normal vector of $C_1$. They showed that $C$ is Mannheim partner curve of $C_1$ if and only if the following equality holds

$$\frac{d\tau}{ds} = \frac{\kappa}{\lambda}\left(1 + \lambda^2\tau^2\right),$$

where $\kappa$ and $\tau$ are the curvature and the torsion of the curve $C$, respectively, and $\lambda$ is a nonzero constant.

In the Euclidean space $E^3$ an oriented line $L$ can be determined by a point $p \in L$ and a normalized direction vector $\vec{a}$ of $L$, i.e. $\|\vec{a}\| = 1$. To obtain the components for $L$, one forms the moment vector $\vec{a}^* = \vec{p} \times \vec{a}$ with respect to the origin point in $E^3$. The two vectors $\vec{a}$ and $\vec{a}^*$ are not independent of one another; they satisfy the relationships $\langle \vec{a}, \vec{a} \rangle = 1$, $\langle \vec{a}, \vec{a}^* \rangle = 0$. The pair $(\vec{a}, \vec{a}^*)$ of the vectors $\vec{a}$ and $\vec{a}^*$, which satisfies those relationships, is called dual unit vector [10]. The most important properties of real vector analysis are valid for the dual vectors. According to the above statement since each dual unit vector corresponds to a line of $E^3$, there is a one-to-one correspondence between the points of a dual unit sphere $\tilde{S}^2$ and the oriented lines of $E^3$. This correspondence is known as E. Study Mapping. As a sequence of that a differentiable curve on dual unit sphere in dual space $\mathbb{D}^3$ represents a ruled surface which is a surface generated by moving of a line $L$ along a curve $\alpha(s)$ in $E^3$ and has the parameterization $r(s,u) = \alpha(s) + uI(s)$, where $\alpha(s)$ is called generating curve and $I$, the direction of the line $L$, is called rulling.

In this paper we study dual Mannheim partner curves in dual space $\mathbb{D}^3$. We will obtain the relationships between the curvatures and the torsions of the dual Mannheim partner curves with respect to each other.

## 2- Preliminaries

Dual numbers had been introduced by W. K. Clifford (1845-1879). A dual number has the form a $\tilde{a} = (a, a^*) = a + \varepsilon a^*$ where $a$ and $a^*$ are real numbers and $\varepsilon = (0,1)$ is called dual unit with $\varepsilon^2 = 0$ where the product of dual numbers a $\tilde{a} = (a, a^*) = a + \varepsilon a^*$, $\tilde{b} = (b, b^*) = b + \varepsilon b^*$ is given by

$$\tilde{a}\tilde{b} = (a, a^*)(b, b^*) = (ab, ab^* + a^*b) = ab + \varepsilon(ab^* + a^*b).$$

We denote the set of dual numbers by $\mathbb{D}$:

$$\mathbb{D} = \{\tilde{a} = a + \varepsilon a^* : a, a^* \in \mathbb{R}, \varepsilon^2 = 0\}.$$

Clifford showed that dual numbers form an algebra, not a field. The pure dual numbers $\varepsilon a^*$ are zero divisors, $(\varepsilon a^*)(\varepsilon b^*) = 0$. However, the other laws of the algebra of dual numbers are the same as the laws of algebra of complex numbers. This means dual numbers form a ring over the real number field.

Now let $f$ be a differentiable dual function. Thus Taylor expansion of dual function $f$ is

$$f(x + \varepsilon x^*) = f(x) + \varepsilon x^* f'(x)$$

where $f'(x)$ is the derivation of $f$ ,. Then we have

$$\sin(x + \varepsilon x^*) = \sin(x) + \varepsilon x^* \cos(x)$$
$$\cos(x + \varepsilon x^*) = \cos(x) - \varepsilon x^* \sin(x) \quad .$$

Let $\mathbb{D}^3$ be the set of all triples of dual numbers, i.e.,

$$\mathbb{D}^3 = \{\tilde{a} = (\tilde{a}_1, \tilde{a}_2, \tilde{a}_3) | \tilde{a}_i \in \mathbb{D}, 1 \le i \le 3\}$$

is a module over the ring $\mathbb{D}$ which is called a $\mathbb{D}$-module or dual space and denoted by $\mathbb{D}^3$. The elements of $\mathbb{D}^3$ are called dual vectors. Thus a dual vector $\vec{\tilde{a}}$ can be written

$$\vec{\tilde{a}} = \vec{a} + \varepsilon \vec{a}^*$$

where $\vec{a}$ and $\vec{a}^*$ are real vectors in $\mathbb{R}^3$. Then for any vectors $\vec{\tilde{a}}$ and $\vec{\tilde{b}}$ in $\mathbb{D}^3$, the scalar product and the vector product are defined by

$$\langle \vec{\tilde{a}}, \vec{\tilde{b}} \rangle = \langle \vec{a}, \vec{b} \rangle + \varepsilon \left( \langle \vec{a}, \vec{b}^* \rangle + \langle \vec{a}^*, \vec{b} \rangle \right)$$

and

$$\vec{\tilde{a}} \times \vec{\tilde{b}} = \vec{a} \times \vec{b} + \varepsilon \left( \vec{a} \times \vec{b}^* + \vec{a}^* \times \vec{b} \right),$$

respectively.

The norm of a dual vector $\vec{\tilde{a}}$ is defined to be

$$\|\vec{\tilde{a}}\| = \|\vec{a}\| + \varepsilon \frac{\langle \vec{a}, \vec{a}^* \rangle}{\|\vec{a}\|} \quad , \quad \vec{a} \neq \vec{0}.$$

A dual vector $\vec{\tilde{a}}$ with norm 1 is called dual unit vector. The set of dual unit vectors is

$$\tilde{S}^2 = \left\{ \vec{\tilde{a}} = \vec{a} + \varepsilon \vec{a}^* \in \mathbb{D}^3 \mid \|\vec{\tilde{a}}\| = 1, \ \vec{a}, \vec{a}^* \in \mathbb{R}^3 \right\}$$

which is called dual unit sphere. Similarly, the set of arbitrary dual unit vectors is

$$\tilde{S}^2(r) = \left\{ \vec{\tilde{a}} = \vec{a} + \varepsilon \vec{a}^* \in \mathbb{D}^3 \mid \langle \vec{\tilde{a}} - \tilde{c}, \vec{\tilde{a}} - \tilde{c} \rangle = \tilde{r}^2 \right\}$$

is a dual sphere with radius $\tilde{r}$ and center $\tilde{c}$.

E. Study used dual numbers and dual vectors in his research on the geometry of lines and kinematics. He devoted special attention to the representation of directed lines by dual unit vectors and defined the mapping that is known by his name: There exists one-to-one correspondence between the vectors of dual unit sphere $\tilde{S}^2$ and the directed lines of space of lines $\mathbb{R}^3$. By the aid of this correspondence, the properties of the spatial motion of a line can be derived. Hence, the geometry of ruled surface is represented by the geometry of dual curves on the dual unit sphere in $\mathbb{D}^3$.

The angle $\tilde{\theta} = \theta + \varepsilon \theta^*$ between two dual unit vectors $\vec{\tilde{a}}$, $\vec{\tilde{b}}$ is called dual angle and defined by

$$\langle \vec{\tilde{a}}, \vec{\tilde{b}} \rangle = \cos(\tilde{\theta}) = \cos(\theta) - \varepsilon \theta^* \sin(\theta).$$

By considering The E. Study Mapping the geometric interpretation of dual angle is that $\theta$ is the real angle between lines $L_1$, $L_2$ corresponding to the dual unit vectors $\vec{\tilde{a}}$, $\vec{\tilde{b}}$, respectively, and $\theta^*$ is the shortest distance between those lines (see figure 1).

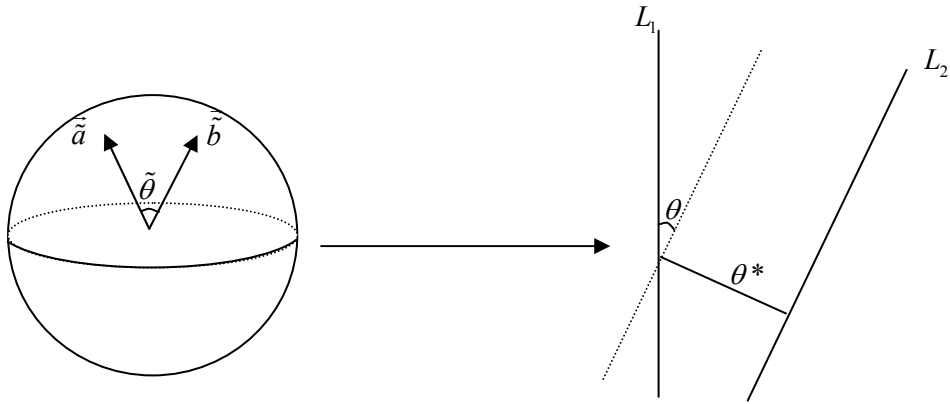

Fig. 1. Dual angle between dual unit vectors and corresponding lines to those vectors

Let $\alpha(t) = (\alpha_1(t), \alpha_2(t), \alpha_3(t))$ and $\alpha^*(t) = (\alpha_1^*(t), \alpha_2^*(t), \alpha_3^*(t))$ be real valued curves in the space $E^3$. Then $\tilde{\alpha}(t) = \alpha(t) + \varepsilon \alpha^*(t)$ is a curve in dual space $\mathbb{D}^3$ and is called dual space curve. If the real valued functions $\alpha_i(t)$ and $\alpha_i^*(t)$ are differentiable then the dual space curve

$$\tilde{\alpha}: I \subset \mathbb{R} \to \mathbb{D}^3$$
$$t \to \tilde{\alpha}(t) = (\alpha_1(t) + \varepsilon \alpha_1^*(t), \alpha_2(t) + \varepsilon \alpha_2^*(t), \alpha_3(t) + \varepsilon \alpha_3^*(t))$$
$$= \alpha(t) + \varepsilon \alpha^*(t)$$

is differentiable in dual space $\mathbb{D}^3$. The real part $\alpha(t)$ of the dual space curve $\tilde{\alpha} = \tilde{\alpha}(t)$ is called indicatrix. The dual arc-length of the dual space curve $\tilde{\alpha}(t)$ from $t_1$ to $t$ is defined by

$$\tilde{s} = \int_{t_1}^{t} \|\vec{\tilde{\alpha}}'(t)\| dt = \int_{t_1}^{t} \|\vec{\alpha}'(t)\| dt + \varepsilon \int_{t_1}^{t} \langle \vec{t}, (\vec{\alpha}^*(t))' \rangle dt = s + \varepsilon s^*,$$

where $\vec{t}$ is unit tangent vector of the indicatrix $\alpha(t)$ which is a real space curve in $E^3$ [11]. From now on we will take the arc length $s$ of $\vec{\alpha}(t)$ as the parameter instead of $t$.

Denote by $\{\vec{\tilde{t}}, \vec{\tilde{n}}, \vec{\tilde{b}}\}$ the moving dual Frenet frame along the dual space curve $\tilde{\alpha}(s)$ in the dual space $\mathbb{D}^3$. Then $\vec{\tilde{t}}, \vec{\tilde{n}}, \vec{\tilde{b}}$ are the dual tangent, the dual principal normal and the dual binormal vector fields, respectively. Then for the curve $\tilde{\alpha}$ the Frenet formulae are given by

$$\frac{d}{d\tilde{s}} \begin{bmatrix} \vec{\tilde{t}} \\ \vec{\tilde{n}} \\ \vec{\tilde{b}} \end{bmatrix} = \begin{bmatrix} 0 & \tilde{\kappa} & 0 \\ -\tilde{\kappa} & 0 & \tilde{\tau} \\ 0 & -\tilde{\tau} & 0 \end{bmatrix} \begin{bmatrix} \vec{\tilde{t}} \\ \vec{\tilde{n}} \\ \vec{\tilde{b}} \end{bmatrix} \qquad (1)$$

the functions $\tilde{\kappa} = \kappa + \varepsilon \kappa^*$ and $\tilde{\tau} = \tau + \varepsilon \tau^*$ are called dual curvature and dual torsion of $\tilde{\alpha}$, respectively. The formulae (1) are called the Frenet formulae of dual curve in $\mathbb{D}^3$ [11].

**3-) Dual Mannheim partner curves**

In this section, we define dual Mannheim partner curves in dual space $\mathbb{D}^3$ and we give some characterization for dual Mannheim partner curves in the same space. Using these relationships, we will comment again Bertrand's, Schell's and Mannheim's theorems.

**Definition 1.** Let $\mathbb{D}^3$ be the dual space with the standard inner product $\langle , \rangle$. If there exists a corresponding relationship between the dual space curves $\tilde{C}$ and $\tilde{C}_1$ such that, at the corresponding points of the dual curves, the principal normal lines of $\tilde{C}$ coincides with the binormal lines of $\tilde{C}_1$, then $\tilde{C}$ is called a dual Mannheim curve, and $\tilde{C}_1$ a dual Mannheim partner curve of $\tilde{C}$. The pair $\{\tilde{C}, \tilde{C}_1\}$ is said to be a dual Mannheim pair (see figure 2).

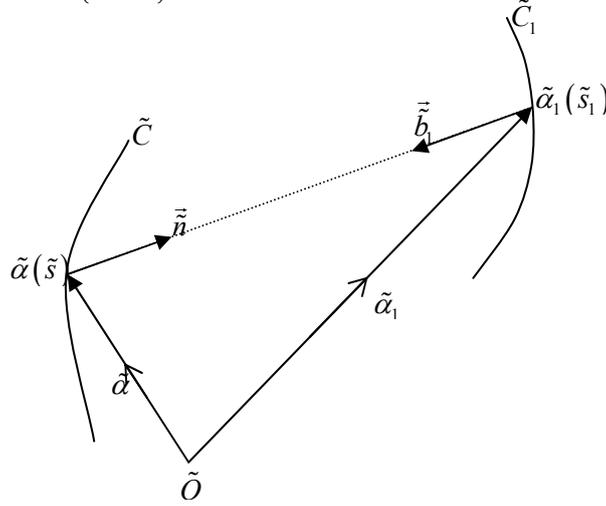

Fig. 2. The Mannheim partner curves

Let $\tilde{\alpha}: \tilde{\alpha}(\tilde{s})$ be a dual Mannheim curve in $\mathbb{D}^3$ parameterized by its arc length $\tilde{s}$ and $\tilde{\alpha}_1: \tilde{\alpha}_1(\tilde{s}_1)$ the dual Mannheim partner curve of with an arc length parameter $\tilde{s}_1$. Denote by $\{\vec{\tilde{t}}(\tilde{s}), \vec{\tilde{n}}(\tilde{s}), \vec{\tilde{b}}(\tilde{s})\}$ the Frenet frame field along $\tilde{\alpha}: \tilde{\alpha}(\tilde{s})$.

In the following theorems, we give a necessary and sufficient condition for a dual space curve to be a Mannheim curve.

**Theorem 1.** A dual space curve in $\mathbb{D}^3$ is a dual Mannheim curve if and only if its curvature $\tilde{\kappa}$ and torsion $\tilde{\tau}$ satisfy the formula $\tilde{\kappa} = \tilde{\lambda}(\tilde{\kappa}^2 + \tilde{\tau}^2)$, where $\tilde{\lambda}$ is never pure dual constant [12].

**Theorem 2.** The distance between corresponding dual points of the dual Mannheim partner curves in $\mathbb{D}^3$ is constant.

**Proof:** From the definition of dual Mannheim curve, we can write

$$\tilde{\alpha}(\tilde{s}) = \tilde{\alpha}_1(\tilde{s}_1) + \tilde{\lambda}(\tilde{s}_1)\vec{\tilde{b}}_1(\tilde{s}_1). \qquad (2)$$

By taking the derivative of this equation with respect to $\tilde{s}_1$ and applying the Frenet formulas, we get

$$\vec{\tilde{t}}\frac{d\tilde{s}}{d\tilde{s}_1} = \vec{\tilde{t}}_1 + \tilde{\lambda}'\vec{\tilde{b}}_1 - \tilde{\lambda}\tilde{\tau}_1\vec{\tilde{n}}_1, \qquad (3)$$

where we use the superscript $(\,'\,)$ to denote the derivative with respect to the arc length parameter of a dual curve. Since the dual vectors $\vec{\tilde{n}}$ and $\vec{\tilde{b}}_1$ are linearly dependent, $\langle \vec{\tilde{t}}, \vec{\tilde{b}}_1 \rangle = 0$. Therefore, from equation (3) we have $\tilde{\lambda}' = 0$. This means that $\tilde{\lambda} = \lambda + \varepsilon \lambda^*$ is a nonzero constant. On the other hand, from the definition of distance function between two points, we get

$$d(\tilde{\alpha}_1(\tilde{s}_1), \tilde{\alpha}(\tilde{s})) = \|\tilde{\alpha}(\tilde{s}) - \tilde{\alpha}_1(\tilde{s}_1)\| = \|\tilde{\lambda}\vec{\tilde{b}}_1\| = |\tilde{\lambda}| = |\lambda| = \text{constant} \neq 0.$$

This completes the proof.

**Theorem 3.** For a dual curve $\tilde{C}$ in $\mathbb{D}^3$, there is a dual curve $\tilde{C}_1$ so that $\{\tilde{C}, \tilde{C}_1\}$ is a dual Mannheim pair.

**Proof:** Since the dual vectors $\vec{\tilde{n}}$ and $\vec{\tilde{b}}_1$ are linearly dependent, equation (2) can be rewritten as

$$\tilde{\alpha}_1(\tilde{s}_1) = \tilde{\alpha}(\tilde{s}) - \tilde{\lambda}(\tilde{s}_1)\vec{\tilde{n}}(\tilde{s}_1). \qquad (4)$$

Since $\tilde{\lambda}$ is a nonzero constant, there is a dual curve $\tilde{C}_1$ for all values of $\tilde{\lambda}$.

Now, we can give the following theorem related to curvature and torsion of the dual Mannheim partner curves.

**Theorem 4.** Let $\{\tilde{C}, \tilde{C}_1\}$ be a dual Mannheim pair in $\mathbb{D}^3$. Then, the torsion of the dual curve $\tilde{C}_1$ is $\tilde{\tau}_1 = \dfrac{\tilde{\kappa}}{\tilde{\lambda}\tilde{\tau}}$.

**Proof:** By considering $\tilde{\lambda}$ is nonzero constant in equation (3), we get

$$\vec{\tilde{t}} = \frac{d\tilde{s}_1}{d\tilde{s}}\vec{\tilde{t}}_1 - \tilde{\lambda}\tilde{\tau}_1 \frac{d\tilde{s}_1}{d\tilde{s}}\vec{\tilde{n}}_1 \; . \qquad (5)$$

By taking the derivative of equation (4) with respect to $\tilde{s}$ and applying the Frenet formulas, we obtain

$$\vec{\tilde{t}}_1 = (1 + \tilde{\lambda}\tilde{\kappa})\frac{d\tilde{s}}{d\tilde{s}_1}\vec{\tilde{t}} - \tilde{\lambda}\tilde{\tau}_1 \frac{d\tilde{s}}{d\tilde{s}_1}\vec{\tilde{b}} \; . \qquad (6)$$

On the other hand, since $\vec{\tilde{n}}$ is coincident with $\vec{\tilde{b}}_1$ in direction, we can write

$$\begin{aligned}\vec{\tilde{t}} &= \cos\tilde{\theta}\,\vec{\tilde{t}}_1 + \sin\tilde{\theta}\,\vec{\tilde{n}}_1 \\ \vec{\tilde{b}} &= -\sin\tilde{\theta}\,\vec{\tilde{t}}_1 + \cos\tilde{\theta}\,\vec{\tilde{n}}_1\end{aligned} \qquad (7)$$

or

$$\vec{\tilde{t}}_1 = \cos\tilde{\theta}\vec{\tilde{t}} - \sin\tilde{\theta}\vec{\tilde{b}}$$
$$\vec{\tilde{n}}_1 = \sin\tilde{\theta}\vec{\tilde{t}} + \cos\tilde{\theta}\vec{\tilde{b}} \;,$$
(8)

where $\tilde{\theta}$ is the dual angle between $\vec{\tilde{t}}$ and $\vec{\tilde{t}}_1$ at the corresponding points of the dual curves of $\tilde{C}$ and $\tilde{C}_1$. By taking into consideration equations (5) and (7), we get

$$\cos\tilde{\theta} = \frac{d\tilde{s}_1}{d\tilde{s}} \;,\quad \sin\tilde{\theta} = -\tilde{\lambda}\,\tilde{\tau}_1 \frac{d\tilde{s}_1}{d\tilde{s}}.$$
(9)

Similarly, from equations (6) and (8), we get

$$\cos\tilde{\theta} = (1+\tilde{\lambda}\,\tilde{\kappa})\frac{d\tilde{s}}{d\tilde{s}_1} \;,\quad \sin\tilde{\theta} = \tilde{\lambda}\,\tilde{\tau}\,\frac{d\tilde{s}}{d\tilde{s}_1}.$$
(10)

Therefore, from equations (9) and (10), we obtain

$$\tilde{\tau}_1 = \frac{\tilde{\kappa}}{\tilde{\lambda}\,\tilde{\tau}}.$$

**Theorem 5.** The dual curve $\tilde{\alpha}$ is a straight line in $\mathbb{D}^3$ if and only if $\tilde{\kappa} = 0$.

**Proof:** It is simply seen that $\tilde{\kappa} = \|\vec{\tilde{t}}'\|$ from equation (1). Thus,

$$\tilde{\kappa} = 0 \Leftrightarrow \|\vec{\tilde{t}}'\| = 0 \Leftrightarrow \vec{\tilde{t}}' = \vec{0}.$$

That the dual vector $\vec{\tilde{t}}'$ is a zero vector means that the dual vector $\vec{\tilde{t}}$ is a constant. Taking this constant vector as $\vec{\tilde{x}}$, we get

$$\tilde{\alpha}(\tilde{s}) = \vec{\tilde{x}}\,\tilde{s} + \tilde{y}.$$

The last equation is a dual straight line equation, that is passing point $y$ and parallel to vector $\vec{x}$ in $\mathbb{R}^3$, where $\tilde{y} = y + \varepsilon y^*$ is a dual point and $\vec{\tilde{x}} = \vec{x} + \varepsilon \vec{x}^*$ is a dual vector.

**Theorem 6.** Let $\tilde{\alpha}$ be a unit speed dual curve in $\mathbb{D}^3$ with nonzero constant curvature $\tilde{\lambda}$. Then $\tilde{\alpha}$ is a plane curve if and only if $\tilde{\tau} = 0$.

**Proof:** Suppose $\tilde{\alpha}$ is a plane curve. Then, there exist a point $\tilde{p}$ and a dual constant vector $\vec{\tilde{q}}$ such that $\langle \tilde{\alpha}(\tilde{s}) - \tilde{p}, \vec{\tilde{q}} \rangle = 0$ for all $\tilde{s}$. Hence

$$\langle \tilde{\alpha}'(\tilde{s}), \vec{\tilde{q}} \rangle = \langle \tilde{\alpha}''(\tilde{s}), \vec{\tilde{q}} \rangle = 0 \qquad \text{for all } \tilde{s}.$$

Thus $\vec{\tilde{q}}$ is always orthogonal to $\vec{\tilde{t}} = \tilde{\alpha}'$ and $\vec{\tilde{n}} = \tilde{\alpha}''/\tilde{\kappa}$. Since $\vec{\tilde{b}}$ has unit length, $\vec{\tilde{b}}$ is also orthogonal to $\vec{\tilde{t}}$ and $\vec{\tilde{n}}$, we get $\vec{\tilde{b}} = \pm \vec{\tilde{q}}/\|\vec{\tilde{q}}\|$. Thus $\vec{\tilde{b}}' = \vec{0}$, and taking into consideration equation (1) $\tilde{\tau} = 0$.

Conversely, suppose $\tilde{\tau} = 0$. Thus $\vec{\tilde{b}}' = \vec{0}$; that is, dual vector $\vec{\tilde{b}}$ is constant. We will assert that $\tilde{\alpha}$ lies in the plane through $\tilde{\alpha}(0)$ and $\vec{\tilde{b}}$ is orthogonal this plane. To prove this, we consider the dual-valued function $\tilde{f}(\tilde{s}) = \langle \tilde{\alpha}(\tilde{s}) - \tilde{\alpha}(0), \vec{\tilde{b}} \rangle$ for all $\tilde{s}$. Then, we get

$$\frac{d\tilde{f}}{d\tilde{s}} = \langle \alpha', \vec{\tilde{b}} \rangle = \langle \vec{\tilde{t}}, \vec{\tilde{b}} \rangle = 0.$$

However, obviously, $\tilde{f}(0) = 0$, so $\tilde{f}(\tilde{s})$ is identically zero. Thus, $\langle \tilde{\alpha}(\tilde{s}) - \tilde{\alpha}(0), \vec{\tilde{b}} \rangle = 0$ for all $\tilde{s}$, which shows that $\tilde{\alpha}$ lies entirely in this plane orthogonal to the binormal of $\tilde{\alpha}$. This completes the proof.

By considering Theorem 4, Theorem 5, and Theorem 6 we can give the following results.

**Corollary 1.** If the dual curve $\tilde{\alpha}$ is straight line in $\mathbb{D}^3$, then dual Mannheim partner curve $\tilde{\alpha}_1$ of the dual curve $\tilde{\alpha}(\tilde{s})$ is a plane curve.

**Corollary 2.** Let $\{\tilde{C}, \tilde{C}_1\}$ be a dual Mannheim pair in $\mathbb{D}^3$. Then, the dual product of torsions $\tilde{\tau}$ and $\tilde{\tau}_1$ at the corresponding points of the dual Mannheim partner curves is not constant.

Namely, Schell's theorem is invalid for the dual Mannheim curves.

**Theorem 7.** Let $\{\tilde{C}, \tilde{C}_1\}$ be a dual Mannheim partner curves in $\mathbb{D}^3$. Between the curvature and the torsion of the dual curve $\tilde{C}$, there is the relationship

$$\tilde{\mu}\tilde{\tau} - \tilde{\lambda}\tilde{\kappa} = 1,$$

where $\tilde{\lambda}$ and $\tilde{\mu}$ are nonzero dual numbers.

**Proof:** From equation (10), we obtain

$$\frac{\cos\tilde{\theta}}{(1+\tilde{\lambda}\tilde{\kappa})} = \frac{\sin\tilde{\theta}}{\tilde{\lambda}\tilde{\tau}}.$$

Arranging this equation, we get

$$\tilde{\lambda}\cot\tilde{\theta}\,\tilde{\tau} - \tilde{\lambda}\tilde{\kappa} = 1,$$

and if we have chose $\tilde{\lambda}\cot\tilde{\theta} = \tilde{\mu}$ for brevity, we see that

$$\tilde{\mu}\tilde{\tau} - \tilde{\lambda}\tilde{\kappa} = 1.$$

**Corollary 3.** Let $\{\tilde{C}, \tilde{C}_1\}$ be a dual Mannheim pair in $\mathbb{D}^3$. Then there exists a linear relationship between $\tilde{\kappa}$ and $\tilde{\tau}$ with constant coefficients.

Namely, Bertrand's theorem is valid for the dual Mannheim partner curves.

**Theorem 8.** Let $\{\tilde{C}, \tilde{C}_1\}$ be a dual Mannheim partner curves in $\mathbb{D}^3$. There are the following equations for the curvatures and the torsions of the curves $\tilde{C}$ and $\tilde{C}_1$:

**i-)** $\tilde{\kappa}_1 = -\dfrac{d\tilde{\theta}}{d\tilde{s}_1}$,

**ii-)** $\tilde{\tau}_1 = \sin\tilde{\theta}\, \tilde{\kappa}\, \dfrac{d\tilde{s}}{d\tilde{s}_1} - \cos\tilde{\theta}\, \tilde{\tau}\, \dfrac{d\tilde{s}}{d\tilde{s}_1}$,

**iii-)** $\tilde{\kappa} = \sin\tilde{\theta}\, \tilde{\tau}_1\, \dfrac{d\tilde{s}_1}{d\tilde{s}}$,

**iv-)** $\tilde{\tau} = -\cos\tilde{\theta}\, \tilde{\tau}_1\, \dfrac{d\tilde{s}_1}{d\tilde{s}}$

**Proof: i-)** By considering equation (7) we can easily seen that $\langle \vec{t}, \vec{t}_1 \rangle = \cos\tilde{\theta}$. Then, by taking the derivative of this equality with respect to $\tilde{s}_1$ and by considering equation (1), we have

$$\left\langle \tilde{\kappa}\vec{n}\, \frac{d\tilde{s}}{d\tilde{s}_1}, \vec{t}_1 \right\rangle + \left\langle \tilde{\kappa}_1\, \vec{n}_1, \vec{t} \right\rangle = -\sin\tilde{\theta}\, \frac{d\tilde{\theta}}{d\tilde{s}_1}.$$

Moreover, since $\vec{n}$ and $\vec{b}_1$ are linearly dependent, we get

$$\left\langle \tilde{\kappa}_1\, \vec{n}_1, \vec{t} \right\rangle = -\sin\tilde{\theta}\, \frac{d\tilde{\theta}}{d\tilde{s}_1}$$

and by considering equation (8), we obtain

$$\tilde{\kappa}_1 = -\frac{d\tilde{\theta}}{d\tilde{s}_1}.$$

In a very similar manner, considering $\langle \vec{n}, \vec{n}_1 \rangle = 0$, $\langle \vec{t}, \vec{b}_1 \rangle = 0$ and $\langle \vec{b}, \vec{b}_1 \rangle = 0$, we can build the proofs of (ii), (iii) and (iv) of the theorem, respectively.

By considering (iii) and (iv) we can give the following result.

**Corollary 4.** $\tilde{\kappa}^2 + \tilde{\tau}^2 = \dfrac{d\tilde{s}_1}{d\tilde{s}}\, \tilde{\tau}_1^2$.

**Theorem 9.** Let $\{\tilde{C}, \tilde{C}_1\}$ be a dual Mannheim partner curves in $\mathbb{D}^3$. Moreover, let the dual points $\tilde{\alpha}(\tilde{s})$, $\tilde{\alpha}_1(\tilde{s}_1)$ be two corresponding dual points of $\{\tilde{C}, \tilde{C}_1\}$ and $\tilde{M}$, $\tilde{M}_1$ be the curvature centers at these points, respectively. Then, the ratio

$$\frac{\|\tilde{\alpha}_1(\tilde{s}_1)\tilde{M}\|}{\|\tilde{\alpha}(\tilde{s})\tilde{M}\|} : \frac{\|\tilde{\alpha}_1(\tilde{s}_1)\tilde{M}_1\|}{\|\tilde{\alpha}(\tilde{s})\tilde{M}_1\|}$$

is not constant.

**Proof:** A circle that lies in the dual osculating plane of the point $\tilde{\alpha}(\tilde{s})$ on the dual curve $\tilde{C}$ and that has the centre $\tilde{M} = \tilde{\alpha}(\tilde{s}) + \frac{1}{\tilde{\kappa}}\vec{\tilde{n}}$ lying on the dual principal normal $\vec{\tilde{n}}$ of the point $\tilde{\alpha}(\tilde{s})$ and the radius $\frac{1}{\tilde{\kappa}}$ far from $\tilde{\alpha}(\tilde{s})$, is called dual osculating circle of the dual curve $\tilde{C}$ in the point $\tilde{\alpha}(\tilde{s})$. Similar definition can be given for the dual curve $\tilde{C}_1$ too (see figure 3).

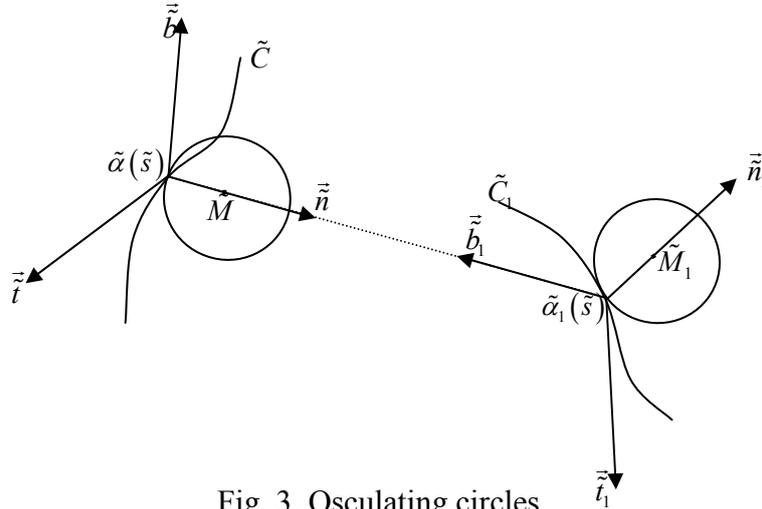

Fig. 3. Osculating circles

Then, we can write

$$\|\tilde{\alpha}(\tilde{s})\tilde{M}\| = \left\|\frac{1}{\tilde{\kappa}}\vec{\tilde{n}}\right\| = \left|\frac{1}{\tilde{\kappa}}\right| = \frac{1}{\kappa} , \qquad \|\tilde{\alpha}_1(\tilde{s}_1)\tilde{M}_1\| = \left\|\frac{1}{\tilde{\kappa}_1}\vec{\tilde{n}}_1\right\| = \left|\frac{1}{\tilde{\kappa}_1}\right| = \frac{1}{\kappa_1}$$

$$\|\tilde{\alpha}(\tilde{s})\tilde{M}_1\| = \left\|\frac{1}{\tilde{\kappa}_1}\vec{\tilde{n}}_1 - \tilde{\lambda}\vec{\tilde{b}}_1\right\| = \frac{\sqrt{1+\tilde{\kappa}_1\tilde{\lambda}^2}}{|\tilde{\kappa}_1|} = \frac{\sqrt{1+\tilde{\kappa}_1\tilde{\lambda}^2}}{\kappa_1} ,$$

$$\|\tilde{\alpha}_1(\tilde{s}_1)\tilde{M}\| = \left\|\left(\frac{1}{\tilde{\kappa}} + \tilde{\lambda}\right)\vec{\tilde{n}}\right\| = \left|\frac{1}{\tilde{\kappa}} + \tilde{\lambda}\right| = \frac{1}{\kappa} + \lambda .$$

Therefore, we obtain

$$\frac{\|\tilde{\alpha}_1(\tilde{s}_1)\tilde{M}\|}{\|\tilde{\alpha}(\tilde{s})\tilde{M}\|} : \frac{\|\tilde{\alpha}_1(\tilde{s}_1)\tilde{M}_1\|}{\|\tilde{\alpha}(\tilde{s})\tilde{M}_1\|} = (1+\lambda\kappa)\sqrt{1+\tilde{\kappa}_1\tilde{\lambda}^2} \neq \text{constant} .$$

Thus, we can give the following corollary.

**Corollary 5.** Mannheim's theorem is invalid for the dual Mannheim partner curve $\{\tilde{C}, \tilde{C}_1\}$ in $\mathbb{D}^3$.

**KAYNAKLAR**